\documentclass[a4,12pt]{amsart}

\usepackage{amssymb}
\usepackage{epsf}
\usepackage{amsmath,epsf}
\usepackage[latin1]{inputenc}
\usepackage{amsfonts}

\numberwithin{equation}{section}
\newtheorem{theo}{Théorème}[section]
\newtheorem{prop}[theo]{Proposition}

\newtheorem{coro}[theo]{Corollaire}

\newtheorem{rema}[theo]{Remarque}
\newtheorem{exem}[theo]{Exemple}

\newenvironment{resume}{\footnotesize\quotation}

\newenvironment{preuve}[0]
    {\par  \nopagebreak[4] \noindent {\em Preuve.}}
    {\nopagebreak[4] \noindent \hfill $\square$ \par \normalsize}

\newcommand{\pref}[1]{(\ref{#1})}

\def\cf{{\it cf. }}
\def\ie{{\it i.e. }}

\def\gx{>_{\rm lex}}

\def\nm{_{n,m}}

\def\C{{\mathbb C}}
\def\N{{\mathbb N}}
\def\Z{{\mathbb Z}}
\def\Q{{\mathbb Q}}

\def\J{{\mathcal J}}

\def\S{{\mathcal S}}

\def\G{{\mathcal G}}

\def\B{{\mathcal B}}

\def\QSym{{QSym}}

\def\H{{\bf H}}
\def\SH{{\bf SH}}

\def\shuffle{{\,\raise
1pt\hbox{$\scriptscriptstyle\cup{\mskip-4mu}\cup$}\,}}

\title[Polynômes super-coinvariants généralisés]{Polynômes quasi-invariants et super-coinvariants pour le groupe symétrique généralisé}

\author{J.-C.~Aval}

\address[Jean-Christophe Aval]{LaBRI\\ Universit\'e Bordeaux 1\\ 351 cours
de 
 la Lib\'eration\\ 33405 Talence cedex\\ France}
\email{aval@labri.fr}
\urladdr{http://dept-info.labri.fr/$\sim$aval}

\date{25 novembre 2003\thanks{Recherche financée par le Programme IHRP de la Commission Européenne, Research Training Network "Algebraic Combinatorics in Europe," grant HPRN-CT-2001-00272.}}

\begin{document} 

\maketitle 

\begin{resume} 
{\sc Résumé.} Un résultat classique d'Artin affirme que l'idéal engendré par les polynô\-mes 
symétriques sans terme constant en $n$ variables est de codimension $n!$. 
L'auteur, F. Bergeron et N. Bergeron ont récemment obtenu un analogue surprenant dans le cas des polynômes
quasi-symétriques. Dans ce cas, l'idéal est de codimension $C_n$, le $n$-ième nombre de Catalan. 
Les polynômes quasi-symétriques peuvent être vus comme invariants d'une action du groupe symétrique $S_n$, définie par F. Hivert.
Le but de ce travail est de généraliser ces travaux au produit en couronne $S_n\wr \Z_m$, connu sous le nom de groupe symétrique généralisé $G\nm$. 
Après avoir défini une action quasi-symétrisante de $G\nm$ sur $\C[x_1,\dots,x_n]$, nous obtenons une description des invariants, et la codimension de l'idéal associé, à savoir $m^n\,C_n$.
\end{resume}

\begin{abstract} 
  A classical result of Artin states that the ideal generated by symmetric polynomials in $n$ variables is of codimension $n!$. 
The author, F. Bergeron and N. Bergeron have recently obtained a surprising analogous in the case of quasi-symmetric polynomials.
In this case, the ideal is of codimension given by $C_n$, the $n$-th Catalan number. 
Quasi-symmetric polynomials are the invariants of a certain action of the symmetric group $S_n$ defined by F. Hivert.
The aim of this work is to generalize these results to the wreath product $S_n\wr \Z_m$, also known as the generalized symmetric group $G\nm$.
We first define a quasi-symmetrizing action of $G\nm$ on $\C[x_1,\dots,x_n]$, then obtain a description of the invariants and the codimension of the associated ideal, which is $m^n\,C_n$.
\end{abstract}

\vskip 0.3cm
\section{Introduction}

Considérons l'alphabet $X$ en $n$ variables $(x_1,\dots,x_n)$. L'espace des polynômes prenant ses variables dans $X$ et à coefficients complexes est noté $\C[X]$. Soit $G\nm=S_n\wr\Z_m$ le produit en couronne du groupe symétriques $S_n$ par le groupe cyclique $\Z_m$. Ce groupe est souvent applelé {\sl groupe symétrique généralisé} (\cf \cite{osima}). On peut se représenter un élément de $G\nm$ comme une matrice carrée de taille $n$ dans laquelle chaque ligne et chaque colonne comporte exactement une entrée non nulle (matrice de pseudo-permutation), et telle que ces entrées non nulles sont des racines $m$-ième de l'unité. L'ordre du groupe $G\nm$ est par conséquent $m^n\,n!$. Si $m=1$, $G\nm$ se réduit au groupe symétrique $S_n$, et si $m=2$, $G\nm$ est le groupe hyperoctaèdral $B_n$, \ie le groupe des permutations signées ou groupe de Weil de type $B$ (voir \cite{lus} pour plus de détails). Le groupe $G\nm$ agit sur les polynômes (action ``classique'') de la façon suivante :
\begin{equation}\label{clas}
\forall g\in G\nm,\ \forall P\in\C[X],\ g.P(X)=P(X.{}^tg),
\end{equation}
o\`u ${}^tg$ est la transposée de la matrice $g$ et $X$ est vu comme vecteur ligne.
Soit
$$Inv\nm=\{P\in\C[X]\ /\ \forall g\in G\nm,\ g.P=P\}$$
l'ensemble des invariants polynomiaux pour cette action $G\nm$. Notons de plus $Inv\nm^+$ l'ensemble des tels polynômes sans terme constant. Nous considérons le produit scalaire suivant sur  $\C[X]$ :
\begin{equation}\label{scal}
\langle P,Q\rangle=P(\partial X)Q(X)\mid_{X=0}
\end{equation}
o\`u $\partial X$ représente $(\partial x_1,\dots,\partial x_n)$ et $X=0$ représente $x_1=\cdots=x_n=0$. Utilisons la notation $\langle S\rangle$ pour l'idéal engendré par une partie $S$ de $\C[X]$. L'espace des polynômes $G\nm$-coinvariants est alors défini par
\begin{eqnarray*}
Cov\nm&=&\{P\in\C[X]\ /\ \forall Q\in Inv\nm^+,\ Q(\partial X)P=0\}\\
&=&\langle Inv\nm^+\rangle^\perp\ \simeq\ \C[X]/\langle Inv\nm^+\rangle.
\end{eqnarray*}
L'égalité et l'isomorphisme précédent ne sont pas triviaux, mais une référence à ce sujet est \cite{orbit}.

Assez de définitions. Un résultat classique d'Artin \cite{artin} affirme que pour $m=1$ (cas du groupe symétrique), la dimension de l'espace coinvariant $\H_n=Cov_{n,1}$ (qualifié dans ce contexte d'harmonique) est égale à $n!$. Chevalley \cite{che} (voir aussi \cite{shetod}) a généralisé ce résultat en montrant que :
\begin{equation}\label{che}
\dim Cov\nm=|G\nm|=m^n\, n!\,.
\end{equation}

Dans le cas du groupe symétrique ($m=1$), d'intéressants résultats ont été obtenus récemment \cite{a9,b1} en étudiant les coinvariants correspondant, non plus aux polynômes symétriques, mais aux polynômes quasi-symétriques. Notre but principal est ici d'obtenir une description analogue dans le cas $m$ quelconque.

L'anneau $QSym$ des fonctions quasi-symétriques a été introduit par Gessel~\cite{ges} dans le cadre du calcul des fonctions génératrices des $P$-partitions~\cite{stanley1}. Ces fonctions quasi-symétriques sont la source de nombreux travaux récents dans plusieurs domaines de la combinatoire \cite{BMSW,MR,NC,stanley2}. 

Dans \cite{a9,b1}, Aval {\it et.~al.} étudient l'espace $\SH_n$ des polynômes super-coinvariants, defini comme l'orthogonal (par rapport à \pref{scal}) de l'idéal engendré par les polynômes quasi-symétriques sans terme constant, et ont prouvé que sa dimension en tant qu'espace vectoriel est donnée par le $n$-ième nombre de Catalan : 
\begin{equation}\label{catsh}
\dim\SH_n=C_n=\frac 1 {n+1} {2n \choose n}.
\end{equation}
Notre principal résultat est une généralisation de l'équation ci-dessus dans le cas du groupe symétrique généralisé $G\nm$.

Cet article est organisé de la façon suivante. Dans la section 2, nous définissons et étudions une action {\em quasi-symétrisante} du groupe $G\nm$ sur $\C[X]$. Les invariants correspondant à cette action sont appelés quasi-invariants et correspondent aux fonctions quasi-symétriques pour $m=1$. La Section 3 est consacrée à la preuve du théorème central (Théorème \ref{main}), qui donne la dimension de l'espace  $SCov\nm$ des polynômes super-coinvariants pour $G\nm$ ; une utilisation (simple) des bases de Gröbner nous permet de calculer une base explicite de $SCov\nm$ (et sa série de Hilbert).


\section{Une action quasi-symétrisante du groupe $G\nm$}

Nous utliserons la notation vectorielle pour les monômes. Plus précisément, pour
$\nu=(\nu_1,\dots,\nu_n)\in\N^n$, nous noterons $X^\nu$ le monôme
\begin{equation}\label{mon}
x_1^{\nu_1} x_2^{\nu_2} \cdots x_n^{\nu_n}.
\end{equation}
Pour tout $P\in\Q[X]$, nous noterons $[X^\nu]\,P(X)$ le coefficient du monôme $X^\nu$ dans $P(X)$.

Notre première tâche est de définir une action quasi-symétrisante du groupe $G\nm$ sur $\C[X]$. Cette action doit répondre aux critères suivants : elle doit se réduire à l'action de Hivert (\cf \cite{hivert}) dans le cas $m=1$ et fournir des polynômes invariants intéressants. Par {\em intéressant}, nous entendons que ces invariants jouissent d'une jolie caractérisation, de même que les polynômes coinvariants. Il s'avère que le choix d'une telle action n'est pas unique. Celle que nous allons étudier est définie de la façon suivante. Soit $A\subset X$ un sous-alphabet de $X$ comportant $l$ variables et $K=(k_1,\dots,k_l)$ un vecteur d'entiers strictement positifs. Nous ordonnerons un vecteur $B$ constitué de variables $x_i$ distinctes multipliées par des racines de l'unité selon l'ordre des variables et le résultat sera noté $(B)_<$.
Voici maintenant comment agit un élément $g\in G\nm$ :
\begin{equation}\label{quasa}
g\bullet A^K=w(g)^{c(K)} {(A.{}^t|g|)_<}^K
\end{equation}
o\`u $w(g)$ est le poids de $g$, \ie le produit de ses entrées non nulles, $|g|$ est la matrice obtenue en prenant le module des entrées de $g$, et le coefficient $c(K)$ est défini ainsi :
$$c(K)=\left\{
\begin{array}{ll}
0&{\rm si \ } \forall i,\ k_i\equiv 0\ [m]\\
1&{\rm sinon.}
\end{array}\right.$$

\begin{exem}\label{goodex}\rm Si $m=3$ et $n=3$, et nous notons $\zeta$ le nombre complexe $\zeta=e^{\frac{2i\pi}{3}}$, alors par exemple
\begin{eqnarray*}
&&\left(
\begin{array}{ccc}
0&0&j\\
1&0&0\\
0&j&0
\end{array}
\right)\,\bullet\,(x_1^2\,x_2)\\
&=&(j^2)^1\left[\left(
\begin{array}{ccc}
0&0&1\\
1&0&0\\
0&1&0
\end{array}
\right)\,.\,(x_1,x_2)\right]_<^{(2,1)}\\
&=&j^2{(x_3,x_1)_<}^{(2,1)}\\
&=&j^2(x_1,x_3)^{(2,1)}\\
&=&j^2\,x_1^2\,x_3.
\end{eqnarray*}
\end{exem}

Un simple calcul permet de vérifier que ceci définit bien une action du groupe $G\nm$ sur $\C[X]$, qui se réduit à l'action de Hivert (\cf \cite{hivert}, Proposition 3.4) dans le cas $m=1$.

Il va de soit que toute définition du coefficient $c(K)$ donne une action du groupe $G\nm$. Parmi cette famille d'actions, celle définie ici respecte les critères énoncés plus haut. Dans le cas particulier du groupe $B_n$, une action apparentée à celle définie ci-dessus (et fournissant les mêmes invariants) est étudiée dans \cite{b5}.

\'Etudions à présent les polynômes invariants et coinvariants de cette action. Nous devons rappeler quelques définitions.

Une {\sl composition} $\alpha =(\alpha_1,\alpha_2, \dots ,\alpha_k)$ de l'entier
positif $d$ est une liste ordonnée d'entiers strictement positifs dont 
la somme vaut $d$. Pour un vecteur $\nu\in\N^n$, notons $c(\nu)$ la composition
obtenue en enlevant les éventuels zéros de $\nu$. Un polynôme $P\in\Q[X]$ est dit
{\sl quasi-symétrique} si et seulement si, pour tous $\nu$ et $\mu$ dans
$\N^n$, nous avons l'égalité
$$[X^\nu]P(X)=[X^\mu]P(X)$$
dès que $c(\nu)=c(\mu)$. L'espace des polynômes quasi-symétriques en $n$
variables est notée $\QSym_n$.

Les polynômes invariants sous l'action \pref{quasa} de $G\nm$ sont qualifiés de {\it quasi-invariant} et l'espace des polynômes quasi-invariants est noté $QInv\nm$, \ie
$$P\in QInv\nm\ \Leftrightarrow\ \forall g\in G\nm,\ g\bullet P=P.$$
Rappelons (\cf \cite{hivert}, Proposition 3.15) que $QInv_{n,1}=QSym_n$. 
La proposition suivante donne une caractérisation élégante de $QInv\nm$.
\begin{prop}\label{propu}
Nous avons
$$P\in QInv_{n,m}\ \Leftrightarrow\ \exists Q\in QSym_n\ /\ P(X)=Q(X^m)$$
o\`u $Q(X^m)=Q(x_1^m,\dots,x_n^m)$.
\end{prop}
\begin{preuve}
Soit $P$ un élément de $QInv\nm$. Notons $\zeta$ la racine $m$-ième de l'unité $\zeta=e^{\frac{2i\pi} m}$ et $g$ l'élément de $G\nm$ dont la matrice est
$$\left(
\begin{array}{cccc}
\zeta&&0&\\
&1&&\\
0&&\ddots&\\
&&&1
\end{array}
\right).$$
Nous observons que l'équation
$$ \frac 1 m (P+g\bullet P+g^2\bullet P+\cdots+g^{m-1}\bullet P)=P$$
implique que tous les exposants apparaissant dans $P$ doivent être des multiples de $m$. Donc il existe un polynôme $Q\in\C[X]$ tel que $P(X)=Q(X^m)$. 
Pour conclure, il suffit de noter que $S_n\subset G\nm$ implique que $P$ est quasi-symétrique, donc $Q$ est aussi quasi-symétrique.

La réciproque est évidente.
\end{preuve}

\vskip 0.3cm

Définissons maintenant les polynômes {\it super-coinvariants} : 
\begin{eqnarray*}
SCov\nm&=&\{P\in\C[X]\ /\ \forall Q\in QInv\nm^+,\ Q(\partial X)P=0\}\\
&=&\langle QInv\nm^+\rangle^\perp\ \simeq\ \C[X]/\langle QInv\nm^+\rangle
\end{eqnarray*}
o\`u le produit scalaire est définie dans \pref{scal}. C'est l'analogue naturel de $Cov_n$ dans le cadre des actions quasi-symétrisantes et $SCov\nm$  se réduit à l'espace des polynômes super-harmoniques $\SH_n$ (\cf \cite{b1}) quand $m=1$. 

\begin{rema}\rm
Il est clair que tout polynôme invariant sous l'action \pref{quasa} est aussi invariant sous \pref{clas}, \ie $Inv\nm\subset QInv\nm$. En prenant l'orthogonal, on obtient $SCov\nm\subset Cov\nm$, ce qui justifie en un sens la terminologie.
\end{rema}

Notre principal résultat est le suivant, qui est une généralisation de \pref{catsh}. On notera aussi une grande similarité avec \pref{che}.

\begin{theo}\label{main}
La dimension de l'espace $SCov\nm$ est donnée par
\begin{equation}\label{maineq}
\dim SCov\nm=m^n\,C_n=m^n\,\frac 1 {n+1} {2n \choose n}.
\end{equation}
\end{theo}

\begin{rema}\rm
Dans le cas du groupe hyperoctaèdral $B_n=G_{n,2}$, C.-O. Chow \cite{chow} a défini une classe $BQSym(x_0,X)$ de polynômes quasi-symétriques du type $B$ en l'alphabet $(x_0,X)$. Il est intéressant de comparer son approche à la notre et de constater qu'elle est bien différente. On observe en particulier que :
$$BQSym(x_0,X)=QSym(X)+QSym(x_0,X).$$
Il est alors assez simple de voir que le quotient $\C[x_0,X]/\langle BQSym^+\rangle$ est isomorphe au quotient $\C[X]/\langle QSym^+\rangle$ étudié dans \cite{b1}. Expliquons cela en quelques mots pour le lecteur intéressé : si $\G$ est la base de Gröbner de $\langle QSym^+\rangle$ construite dans \cite{b1} (on pourra aussi se reporter à la section suivante), alors l'ensemble $\{x_0,\G\}$ est une base de Gröbner de $\langle BQSym^+\rangle$ (toute syzygie est réductible en vertu du premier principe de Buchberger, \cf \cite{CLO}).
\end{rema}

La section suivante est consacrée à la preuve du Théorème \ref{main}.


\section{Preuve du théorème principal}

Nous allons prouver le Théorème \ref{main} en construisant une base explicite du quotient $\C[X]/\langle QInv\nm^+\rangle$. Ceci repose en grande partie sur les travaux \cite{a9,b1}. Nous allons cependant rappeler ici brièvement les éléments nécessaires.

Commençons par rappeler (\cf \cite{b1}) la bijection suivante qui associe à tout vecteur $\nu\in\N^n$ un chemin $\pi(\nu)$ dans le plan $\N\times\N$. Ce chemin fait des pas Nord ou Est et est défini ainsi : si
$\nu=(\nu_1,\dots,\nu_n)$, le chemin $\pi(\nu)$ est
$$(0,0)\rightarrow(\nu_1,0)\rightarrow(\nu_1,1)\rightarrow(\nu_1+\nu_2,1)
\rightarrow(\nu_1+\nu_2,2)\rightarrow\cdots\ \ \ \ \ \ \ \ \ \ \ \ \ \ \ \ \ \
\ \ \ $$
\vskip -0.6cm
$$\ \ \ \ \ \ \ \ \ \ \ \ \ \ \ \ \ \ \ \ \ \ \ \ \ \ \ \ \ \ \ \ \ \ \ \ \
\rightarrow(\nu_1+\cdots+\nu_n,n-1)\rightarrow(\nu_1+\cdots+\nu_n,n).$$
Par exemple le chemin associé à  $\nu=(2,1,0,3,0,1)$ est donnée à la Figure 1.

\vskip 0.2cm
\begin{figure}[htbp]
\epsfxsize=7cm
$$\epsfbox{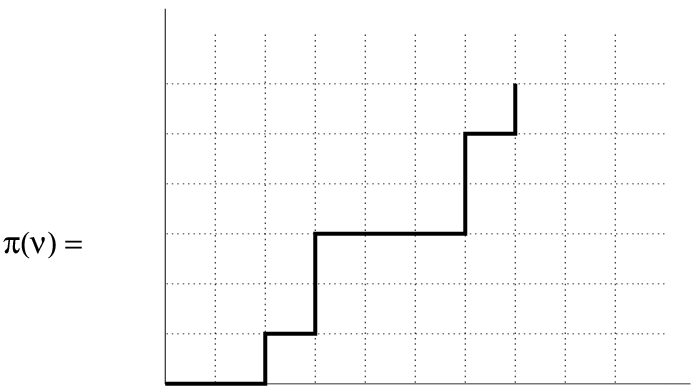}$$
\caption{}
\end{figure}


\vskip 0.2cm

Nous distinguons alors deux types de chemins, suivant leur comportement par rapport à la digonale $y=x$.
Si le chemin reste au-dessus de la diagonale, nous l'appelons un {\em chemin de Dyck}, et
qualifions le vecteur correspondant de {\sl Dyck}. Sinon, nous disons que le chemin
(et le vecteur associé) est {\sl transdiagonal}. Par exemple
$\eta=(0,0,1,2,0,1)$ est de Dyck et $\varepsilon=(0,3,1,1,0,2)$ est
transdiagonal (\cf Figure 2).

\vskip 0.2cm
\begin{figure}[htbp]
\epsfxsize=10cm
$$\epsfbox{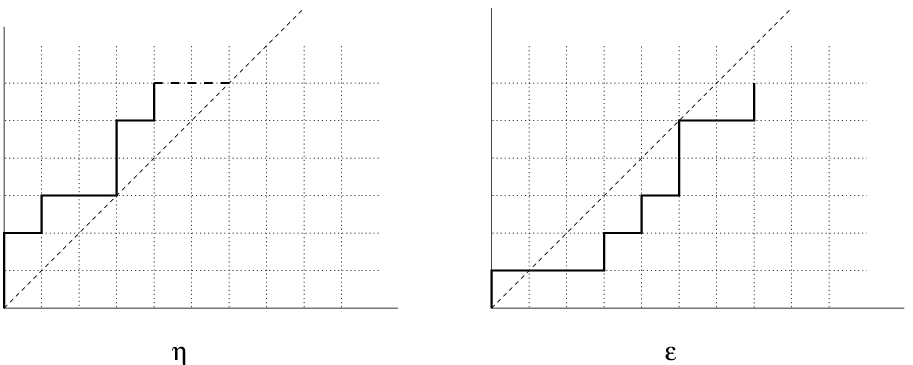}$$
\caption{}
\end{figure}


\vskip 0.2cm

Nous avons alors le résultat suivant qui généralise le Theorem 4.1 de \cite{b1} et qui implique de façon évidente le Théorème \ref{main}.
\begin{theo}\label{main2}
L'ensemble de monômes
$$\B\nm=\{(X_n)^{m \,\eta+\alpha} /\ \pi(\eta)\ {\it est\ un\ chemin\ de\ Dyck},\ 0\le\alpha_i<m\}$$
est une base du quotient $\C[X_n]/\langle QInv_{n,m}^+\rangle$.
\end{theo}

Pour prouver ce résultat, le but est de construire une base de Gröbner pour l'idéal $\J\nm=\langle QInv_{n,m}^+\rangle$. Nous utiliserons les résultats de \cite{a9,b1}.

L'{\it ordre lexicographique} sur les monômes est défini par
\begin{equation}
X^\nu\gx X^\mu\qquad{\rm ssi}\qquad \nu\gx \mu,
\end{equation}
si et seulement si la première composante non nulle du vecteur
$\nu-\mu$ est positive.

Pour toute partie $\S$ de $\Q[X]$ et tout entier strictement positif $m$, introduisons $\S^m=\{P(X^m)\ ,\ P\in\S\}$. Nous noterons $G(I)$ l'unique base de Gröbner réduite (\cf \cite{CLO}) d'un idéal $I$. 
Le lemme suivant, en dépit d'une preuve très simple, est non seulement crucial dans notre contexte mais peut également fournir une preuve instantanée de l'égalité \pref{che}.

\begin{prop}\label{ppp}
Avec les notations précédentes,
\begin{equation}
G(\langle\S^m\rangle)=G(\langle\S\rangle)^m.
\end{equation}
\end{prop} 
\begin{preuve}
C'est une application directe du critère de Buchberger (\cf \cite{CLO} pour une présentation claire du sujet). En effet, si pour toute paire $g,\,g'$ d'éléments de $G(\langle S\rangle)$, la syzygie
$$S(g,g')$$
se réduit à zéro (par hypothèse $G(\langle S\rangle)$ est une base de Gröbner), alors la syzygie
$$S(g(X^m),g'(X^m))$$
se réduit aussi à zéro dans $G(\langle S^m\rangle)$ par exactement le même calcul.
\end{preuve}

\vskip 0.3cm

Rappelons que dans \cite{a9} est construite une famille $\G$ de polynômes $G_\varepsilon$, 
indexés par les vecteurs transdiagonaux $\varepsilon$, possédant les caractéristiques suivantes :
\begin{itemize}
\item le monôme dominant de $G_\varepsilon$ est $LM(G_\varepsilon)=X^\varepsilon$ ;
\item $\G$ est une base de Gröbner de $\J_{n,1}$.
\end{itemize}

Le résultat suivant est alors une conséquence des Propositions \ref{propu} et \ref{ppp}. 
\begin{prop}
L'ensemble $\G^m$ est une base de Gröbner de l'idéal $\J\nm$.
\end{prop}

Pour conclure la preuve du Théorème \ref{main2}, il est suffisant de remarquer que l'ensemble des monômes non divisibles par un monôme dominant d'un élément de $\G^m$, \ie par un $X^{m\varepsilon}$ pour $\varepsilon$ transdiagonal, sont précisément les monômes apparaissant dans $\B\nm$.

Comme corollaire du Théorème \ref{main2}, nous obtenons une formule explicite pour la série de Hilbert de $SCov\nm$. Pour $k\in\N$, notons $SCov\nm^{(k)}$ la projection
\begin{equation}
SCov\nm^{(k)}=SCov\nm\,\cap\,\Q^{(k)}[X]
\end{equation}
o\`u $\Q^{(k)}[X]$ est l'espace des polynômes homogènes de degré $k$ (incluant le polynôme nul). 

Notons $F\nm(t)$ la série de Hilbert de $SCov\nm$, \ie
\begin{equation}
F\nm(t)=\sum_{k\ge 0} \dim SCov\nm^{(k)}\, t^k.
\end{equation}

Rappelons que dans \cite{b1} est calculée la série de Hilbert $F_{n,1}$ :
\begin{equation}\label{ff}
F_{n,1}(t)=F_n(t)=\sum_{k=0}^{n-1}\frac{n-k}{n+k}{n+k
\choose k} t^k
\end{equation}
o\`u apparaît le nombre de chemins de Dyck ayant un nombre fixé (à savoir $n-k$) de facteurs (\cf \cite{Kreweras}).

Le Théorème \ref{main2} implique alors le

\begin{coro}\label{gc}
Avec les notations de \pref{ff}, la série de Hilbert de $SCov\nm$ est donnée par
$$F\nm(t)=\frac{1-t^m}{1-t}F_n(t^m)=\frac{1-t^m}{1-t}\sum_{k=0}^{n-1}\frac{n-k}{n+k}{n+k
\choose k} t^{mk}.$$
De cette formule, on peut déduire la formule close suivante
$$\sum_n F_{n,m}(t)\,x^n=\frac{(1\!-\!t)-\sqrt{(1\!-\!t)(1-t-4t^mx(1-t^m))}-2x(1\!-\!t^m)}{(1-t)(2t^m-1)-x(1-t^m)}\raise 2pt\hbox{.}$$
\end{coro}

\vskip 0.7cm

\section*{Extended abstract in English}

In this extended abstract, the equations, propositions, $\dots$ are numbered as in the French part.

\vskip 0.2cm

Let $X$ denote the alphabet in $n$ variables $(x_1,\dots,x_n)$ and $\C[X]$ denote the space of polynomials with complex coefficients in the alphabet $X$. Let $G\nm=S_n\wr\Z_m$ denote the wreath product of the symmetric group $S_n$ by the cyclic group $\Z_m$, sometimes known as the {\sl generalized symmetric group} (\cf \cite{osima}). It may be seen as the group of $n\times n$ matrices in which each row and each column has exactly one non-zero entry (pseudo-permutation matrices), and such that the non-zero entries are $m$-th roots of unity. The order of $G\nm$ is $m^n\,n!$. When $m=1$, $G\nm$ reduces to the symmetric group $S_n$, and when $m=2$, $G\nm$ is the hyperoctahedral group $B_n$, \ie the group of signed permutations, which is the Weyl group of type $B$. The group $G\nm$ acts classically on $\C[X]$ by the rule
$$\forall g\in G\nm,\ \forall P\in\C[X],\ g.P(X)=P(X.{}^tg),\leqno{(1.1)}$$
where ${}^tg$ is the transpose of the matrix $g$ and $X$ is considered as a row vector.
Let
$$Inv\nm=\{P\in\C[X]\ /\ \forall g\in G\nm,\ g.P=P\}$$
denote the set of $G\nm$-invariant polynomials. Let us denote by $Inv\nm^+$ the set of such polynomials with no constant term. We consider the following scalar product on $\C[X]$:
$$\langle P,Q\rangle=P(\partial X)Q(X)\mid_{X=0}\leqno{(1.2)}$$
where $\partial X$ stands for $(\partial x_1,\dots,\partial x_n)$ and $X=0$ stands for $x_1=\cdots=x_n=0$. Let $\langle S\rangle$ denote the ideal generated by a subset $S$ of $\C[X]$. The space of $G\nm$-coinvariant polynomials is then defined by
\begin{eqnarray*}
Cov\nm&=&\{P\in\C[X]\ /\ \forall Q\in Inv\nm^+,\ Q(\partial X)P=0\}\\
&=&\langle Inv\nm^+\rangle^\perp\ \simeq\ \C[X]/\langle Inv\nm^+\rangle.
\end{eqnarray*}
The previous equality and isomorphism are not obvious, and a reference on that topic is \cite{orbit}.

A classical result of Chevalley \cite{che} (see also \cite{shetod}) states the following equality:
$$\dim Cov\nm=|G\nm|=m^n\, n!\leqno{(1.3)}$$
which reduces when $m=1$ to the theorem of Artin \cite{artin} that the dimension of the harmonic space $\H_n=Cov_{n,1}$ (\cf \cite{orbit}) is $n!$.

Our aim is to give an analogous result in the case of a quasi-symmetrizing action. The ring $QSym$ of quasi-symmetric functions was introduced by 
 Gessel~\cite{ges} as a source 
 of generating functions for $P$-partitions~\cite{stanley1} and  
 appears in more and more combinatorial 
   contexts~\cite{BMSW,stanley2}.  

In \cite{a9,b1}, Aval {\it et.~al.} investigated the space $\SH_n$ of super-coinvariant polynomials for the symmetric group, defined as the orthogonal (with respect to \pref{scal}) of the ideal generated by quasi-symmetric polynomials with no constant term, and proved that its dimension as a vector space equals the $n$-th Catalan number: 
$$\dim\SH_n=C_n=\frac 1 {n+1} {2n \choose n}.\leqno{(1.4)}$$
Our main result is a generalization of the previous equation in the case of super-coinvariant polynomials for the group $G\nm$.

In Section 2, we define and study a {\em quasi-symmetrizing} action of $G\nm$ on $\C[X]$. 
We want this action to give Hivert's action (\cf \cite{hivert}) in the case $m=1$ and to give interesting invariants and coinvariants. Such an action is not unique and we study the one defined as follows. Let $A\subset X$ be a subset of $X$ with $l$ variables and $K=(k_1,\dots,k_l)$ a vector of positive ($>0$) integers. We order a vector $B$ consisting of distinct variables $x_i$ multiplied by roots of unity with respect to the variable order and the result is denoted by $(B)_<$.
Now the quasi-symmetrizing action of $g\in G\nm$ is given by (see also Example \ref{goodex}):
$$g\bullet A^K=w(g)^{c(K)} {(A.{}^t|g|)_<}^K\leqno{(2.2)}$$
where $w(g)$ is the weight of $g$, \ie the product of its non-zeero entries, $|g|$ is the matrix obtained by taking the module of the entries of $g$, and the coefficient $c(K)$ is defined as:
$$c(K)=\left\{
\begin{array}{ll}
0&{\rm if \ } \forall i,\ k_i\equiv 0\ [m]\\
1&{\rm if\ not.}
\end{array}\right.$$

We now study invariants and coinvariants relative to this action.

A {\sl composition} $\alpha =(\alpha_1,\alpha_2, \dots ,\alpha_k)$ of the positive integer $d$ is an ordered list of nonnegative  integers whose sum equals $d$. For $\nu\in\N^n$, let $c(\nu)$ denote the composition
obtained by erasing the zero parts of $\nu$. A polynomial $P\in\Q[X]$ is
{\sl quasi-symmetric} if and only if, for any $\nu,\ \mu\,\in\N^n$, we have
$$[X^\nu]P(X)=[X^\mu]P(X)$$
as soon as $c(\nu)=c(\mu)$. The space of quasi-symmetric polynomials in $n$ variables is denoted by $\QSym_n$.

The invariant of the action \pref{quasa} are said {\it quasi-invariant} and their space is denoted by $QInv\nm$.
Recall (\cf \cite{hivert}, Proposition 3.15) that $QInv_{n,1}=QSym_n$. 
The following proposition gives a characterization of quasi-invariant polynomials.

\vskip 0.2cm
\noindent
{\bf Proposition \ref{propu}.} 
$$P\in QInv_{n,m}\ \Leftrightarrow\ \exists Q\in QSym_n\ /\ P(X)=Q(X^m) \ \ {\rm with}\ Q(X^m)=Q(x_1^m,\dots,x_n^m).$$

We now define {\it super-coinvariant} polynomials: 
\begin{eqnarray*}
SCov\nm&=&\{P\in\C[X]\ /\ \forall Q\in QInv\nm^+,\ Q(\partial X)P=0\}\\
&=&\langle QInv\nm^+\rangle^\perp\ \simeq\ \C[X]/\langle QInv\nm^+\rangle
\end{eqnarray*}

Our main result is the following, which is a generalization of \pref{catsh}, but also shows similarity to \pref{che}.

\vskip 0.2cm
\noindent
{\bf Theorem \ref{main}.} 
The dimension of $SCov\nm$ is given by
$$\dim SCov\nm=m^n\,C_n=m^n\,\frac 1 {n+1} {2n \choose n}.\leqno{(2.3)}$$

\vskip 0.2cm

The Section 3 is devoted to the proof of this result. More precisely, we construct an explicit basis of the quotient $\C[X_n]/\langle QInv_{n,m}^+\rangle$ (\cf Theorem \ref{main2}) from which we deduce its Hilbert series, given by Corollary \ref{gc}.

To do this, we use the results of \cite{a9,b1} to construct a Gröbner basis of the ideal $\langle QInv_{n,m}^+\rangle$.
The important point of this proof is the Proposition \ref{ppp}, so we shall say a few words about it.

The {\it  lexicographic order} on monomials is defined by
$$X^\nu\gx X^\mu\qquad{\rm iff}\qquad \nu\gx \mu,\leqno{(3.1)}$$
if and only if the first non-zero entry of $\nu-\mu$ is positive.

For any $\S\subseteq\Q[X]$ and $m\in\N^*$, we introduce $\S^m=\{P(X^m)\ ,\ P\in\S\}$. We denote by $G(I)$ the unique reduced monic Gröbner basis (\cf \cite{CLO}) of an ideal $I$. The following lemma has a simple proof but is the crucial tool in our context. Furthermore, despite its simplicity, it can provide a proof of \pref{che} in a few lines.

\vskip 0.2cm
\noindent
{\bf Proposition \ref{ppp}.}
With the previous notations,
$$G(\langle\S^m\rangle)=G(\langle\S\rangle)^m.\leqno{(3.2)}$$
\proof
This is a direct consequence of Buchberger's criterion. Indeed, if for every pair $g,\,g'$ in $G(\langle S\rangle)$, the syzygy $S(g,g')$
reduces to zero, then the syzygy $S(g(X^m),g'(X^m))$
also reduces to zero in $G(\langle S^m\rangle)$ by exactly the same computation.
\endproof



\vskip 0.3cm

\begin{thebibliography}{10}

\bibitem{artin} 
 {\sc E. Artin}, 
    {\em Galois Theory},
 Notre Dame Mathematical Lecture {\bf 2} (1944), Notre Dame, IN.

\bibitem{a9} 
   {\sc J.-C. Aval and N. Bergeron}, {\em Catalan Paths and Quasi-Symmetric
Functions}, Proc. Amer. Math. Soc. {\bf 131} (2003), 1053-1062.

\bibitem{b1} 
   {\sc J.-C. Aval, F. Bergeron and N. Bergeron}, {\em Ideals of Quasi-Symmetric Functions and
Super-Coinvariant Polynomials for $S_n$}, Adv. in Math., à paraître.

\bibitem{b5} 
   {\sc J.-C. Aval and N. Bergeron}, {\em Quasi-symmetric polynomials and Temperley-Lieb algebra of type B}, en préparation.

\bibitem{Berg2} 
 {\sc F. Bergeron, A. Garsia and G. Tesler}, {\em Multiple Left Regular
 Representations Generated by Alternants}, J. of Comb. Th., Series A, {\bf
91, 
 1-2}  (2000), 49--83.

\bibitem{BMSW} 
 {\sc N.~Bergeron, S.~Mykytiuk, F.~Sottile, and S.~van Willigenburg}, {\em
    Pieri Operations on Posets},
 \newblock J. of Comb. Theory, Series A, {\bf 91} (2000), 84--110 .

\bibitem{che}
{\sc C. Chevalley}, {\em Invariants of finite groups generated by reflections}, Amer. J. Math., {\bf 77} (1955), 778--782.

\bibitem{chow}
{\sc C.-O. Chow}, {\em Noncommutative Symmetric Functions of Type $B$}, Thesis, Massachusetts Institute of Technology, 2001.

\bibitem{CLO} 
 {\sc D. Cox, J. Little and D. O'Shea}, {\em Ideals, Varieties, and
Algorithms}, 
 Springer-Verlag, New-York, 1992.

\bibitem{orbit}
{\sc A. Garsia, M. Haiman}, {\em Orbit Harmonics and Graded Representations},
\'Editions du Lacim, à paraître.

\bibitem{NC} 
   {\sc I.~Gelfand, D.~Krob, A.~Lascoux, B.~Leclerc, V.~Retakh, and J.-Y.
 Thibon}, 
    {\em Noncommutative symmetric functions}, Adv. in Math., {\bf 112}
(1995), 
   ~218--348. 

\bibitem{ges} 
 {\sc I.~Gessel}, {\em Multipartite ${P}$-partitions and products of skew
    {S}chur functions}, in Combinatorics and Algebra (Boulder, Colo., 1983),
    C.~Greene, ed., vol.~34 of Contemp. Math., AMS, 1984, pp.~289--317.

\bibitem{hivert} 
   {\sc F. Hivert}, {\em Hecke algebras, difference operators, and
 quasi-symmetric functions\/}, Adv. in Math., {\bf 155} (2000),
 181--238. 

\bibitem{Kreweras}
{\sc G. Kreweras}, {\em Sur les \'eventails de segments}, Cahiers du BURO,
{\bf 15} (1970), 3--41.

\bibitem{lus}
{\sc G. Lusztig}, {\em Irreducible representationsof finite reflections groups}, Invent. Math., {\bf 43} (1977), 125--175.

\bibitem{mac} 
 {\sc I.~Macdonald}, {\em Symmetric Functions and Hall Polynomials}, Oxford
    Univ.~Press, 1995,
 \newblock second edition.

\bibitem{MR}  {\sc C. Malvenuto and C. Reutenauer}, {\em Duality
    between quasi-symmetric functions
    and the {S}olomon descent algebra},
    Journal of  Algebra, \textbf{177} (1995), 967--982.

\bibitem{osima}
{\sc M. Osima}, {\em On the representations of the generalized symmetric group}, Math. J. Okayama Univ., {\bf 4} (1954), 39--56.

\bibitem{shetod}
{\sc G.C. Shephard and J.A. Todd}, {\em Finite unitary reflection groups}, CAnadian J. Math {\bf 6} (1954), 274--304.

\bibitem{stanley1} 
 {\sc R.~Stanley}, {\em Enumerative Combinatorics, Vol.~1}, Wadsworth and
    Brooks/Cole, 1986.

\bibitem{stanley2} 
   {\sc R.~Stanley}, {\em Enumerative Combinatorics Vol.~2}, no.~62 in
 Cambridge 
    Studies in Advanced Mathematics, Cambridge University Press, 1999.
 \newblock Appendix 1 by Sergey Fomin.


\end{thebibliography}
\end{document}